\documentclass{aims}
\usepackage{amsmath}
  \usepackage{paralist}
  \usepackage{graphics} 
  \usepackage{epsfig} 
 \usepackage[colorlinks=true]{hyperref}
\hypersetup{urlcolor=blue, citecolor=red}

\def\currentvolume{X}
 \def\currentissue{X}
  \def\currentyear{200X}
   \def\currentmonth{XX}
    \def\ppages{X--XX}
     \def\DOI{10.3934/xx.xx.xx.xx}

\newtheorem{theorem}{Theorem}[section]
\newtheorem{corollary}[theorem]{Corollary}

\newtheorem{lemma}[theorem]{Lemma}
\newtheorem{proposition}[theorem]{Proposition}

\theoremstyle{definition}
\newtheorem{definition}[theorem]{Definition}
\newtheorem{remark}[theorem]{Remark}

\usepackage{mathrsfs}
\newtheorem{example}{Example}

\title[Characterization of minimizable functionals] 
      {Characterization of minimizable Lagrangian action functionals and a dual Mather theorem}
\author[Rodolfo R\'ios-Zertuche]{}


 \email{rodolforiosz@gmail.com}

\thanks{The research leading to these results has received funding from the European Research Council
under the European Union's Seventh Framework Programme (FP/2007-2013) / ERC Grant Agreement
307062. The author acknowledges the support of ANR-3IA Artificial and
Natural Intelligence Toulouse Institute.}

\begin{document}
\maketitle

\centerline{\scshape Rodolfo R\'ios-Zertuche}
\medskip
{\footnotesize
   \centerline{Universit\'e Toulouse 1 Capitole}
 \centerline{Unit\'e Mixte de Recherche ``Toulouse School of Economics -- Research''}
   \centerline{ 1 Esplanade de l'Universit\'e, 31080 Toulouse Cedex 08, France}
   \medskip
  \centerline{ \emph{In grateful memory of Professor John N. Mather}}

} 

\bigskip

 \centerline{(Communicated by the associate editor name)}

\begin{abstract}
We show that a necessary and sufficient condition for a smooth function on the tangent bundle of a manifold to be a Lagrangian density whose action can be minimized is, roughly speaking, that it be the sum of a constant, a nonnegative function vanishing on the support of the minimizers, and an exact form. 

We show that this exact form corresponds to the differential of a Lipschitz function on the manifold that is differentiable on the projection of the support of the minimizers, and its derivative there is Lipschitz. 
This function generalizes the notion of subsolution of the Hamilton-Jacobi equation that appears in weak \textsc{kam} theory, and the Lipschitzity result allows for the recovery of Mather's celebrated 1991 result as a special case. We also show that our result is sharp with several examples. 

Finally, we apply the same type of reasoning to an example of a finite horizon Legendre problem in optimal control, and together with the Lipschitzity result we obtain the Hamilton--Jacobi--Bellman equation and the Maximum Principle. 
\end{abstract}


\emph{This is a corrected version of the published manuscript. For details, please refer to the list of changes is presented at the end. This revision was done in November 2021.}

\section{Introduction}
\label{sec:intro}
In this section, we will recall three areas to which our analysis techniques can be applied, and we sketch the respective result we obtain for each of them and its significance.

\subsection{Lagrangian dynamics}

Let $M$ be a smooth manifold, and let $L\colon TM\to{\mathbb{R}}$ be a time-independent Lagrangian density of class $C^2$. When considering the case in which a closed curve $\gamma\colon[0,T]\to{\mathbb{R}}$ minimizes its action 
\begin{equation}\label{eq:actiontoy}
 A_L(\gamma)=\frac1{T}\int_{0}^{T} L(\gamma(t),\gamma'(t))\,dt,
\end{equation}
our main result, Theorem \ref{thm:mather}, reduces to the following statement.
\begin{theorem}[Toy version of Theorem \ref{thm:mather}]\label{thm:toymather}
 Let $L\colon TM\to{\mathbb{R}}$ be a function of class $C^2(TM)$, and let $\gamma\colon TM\to{\mathbb{R}}$ be an absolutely-continuous closed curve, $\gamma(0)=\gamma(T)$. Then $\gamma$ minimizes the action $A_L$ defined in \eqref{eq:actiontoy} if, and only if, there is a Lipschitz function $f\colon M\to{\mathbb{R}}$ such that
 \begin{enumerate}
 \item\label{it:decompositontoy}
  we have
  \begin{equation}\label{eq:decompositiontoy}
   L\geq df+A_L(\gamma)
  \end{equation}
  wherever the differential $df$ is defined, 
 \item\label{it:differentiabletoy}
  $f$ is differentiable on the image $\gamma([0,T])$,
 \item\label{it:lipschitzitytoy} the map $\gamma(t)\mapsto df_{\gamma(t)}$ from the image $\gamma([0,T])$ to the cotangent bundle is Lipschitz on $\gamma([0,T])$, 
 \item\label{it:decompositionequalitytoy}
  \eqref{eq:decompositiontoy} is actually an equality throughout the image of $(\gamma,\gamma')$, in other words,
  \begin{equation}\label{eq:decompositionequalitytoy}
   L(\gamma(t),\gamma'(t))=df_{\gamma(t)}(\gamma'(t))+A_L(\gamma),\quad t\in[0,T],
  \end{equation}
  and
 \item\label{it:tangentialtoy}
  the momenta of $\gamma$ coincide with the differential $df$, that is,
  \[\frac{\partial L}{\partial v}(\gamma(t),\gamma'(t))=df_{\gamma(t)},\quad t\in[0,T],\]
  where $\partial/\partial v$ denotes differentiation in the fiberwise direction of $T_{\gamma(t)}M$.
 \end{enumerate}
\end{theorem}

In Theorem \ref{thm:mather}, the minimization is considered in a more general context: we allow closed measures on $TM$ ---rather than just curves--- to be the candidates for minimization, and we allow them to have non-trivial boundary. 

But the crux of the matter is already present in the toy version, Theorem \ref{thm:toymather}. Those familiar with Mather's theory \cite{matheractionminimizing91} will recognize item \ref{it:lipschitzitytoy} as analogue to Mather's Lipschitz regularity result, except here the theorem gives regularity of the momenta $\partial L/\partial v$; compare with item \ref{it:tangentialtoy}. Those familiar with Fathi--Siconolfi's  weak \textsc{kam} theory \cite{fathibook,fathisiconolfi04,fathisiconolfi05} will recognize properties in items \ref{it:decompositontoy} and \ref{it:decompositionequalitytoy} as meaning that $f$ is a critical subsolution of the Hamilton-Jacobi equation. Those familiar with Ma\~n\'e's theory \cite{contrerasiturriagabook} will recognize the real number $-A_L(\gamma)$ as Ma\~n\'e's critical value, as in this case it coincides with $-\inf_\eta A_L(\eta)$, where the infimum is taken over all absolutely continuous closed curves $\eta$.

The twist here is that we do not require the function $L$ to be either convex, super-linear, bounded from below, Tonelli, quasi-convex, or coercive; our key assumption instead is that its action is minimizable. We also do not assume that the action functional is semi-continuous, or sequentially continuous. Furthermore, we do not assume that the minimizers are invariant under the Euler--Lagrange flow, and this flow may in fact not be well defined. The minimizers may not even enjoy a graph property, that is, a single velocity vector $v$ for each point $x$ in their projected support $\pi(\operatorname{supp}\mu)$. In particular, the minimizers need not define a foliation or a lamination on $M$. 

Similarly, the Hamiltonian flow and the Hamiltonian itself (that is, the Legendre--Fenchel convex conjugate of $L$) may not be well defined. 
However, we do find an energy conservation principle (Corollary \ref{coro:energyconservation}).

Since we do not restrict the set of possible minimizers very much, the regularity part of our result does not occur on them (see Example \ref{ex:exactform}) --- instead, it occurs in the momenta $\partial L/\partial v$. This motivates the the term ``dual'' in the title of the paper, as Mather's original result concerned the velocities of the minimizers and the result presented here concerns the dual object, namely, their momenta $\partial L/\partial v$.

A version of Theorem \ref{thm:mather} for time-dependent Lagrangians is given in Corollary \ref{cor:mathertimedependent}.  Before that, we establish a preliminary result that requires less regularity on the Lagrangian $L$ and gives a weaker characterization result; this is Theorem \ref{thm:characterization}, whose version in the time-dependent setting is Corollary \ref{coro:timedepcoarse}; their proof is the goal of Section \ref{sec:coarsecharact}. %

In the literature, results stating that an object in the tangent bundle $TM$ is contained in the graph of a section of the bundle, or equivalently, that it intersects each fiber at no more than one point, are known as ``graph theorems''; see for example \cite{matheractionminimizing91,carneiroruggiero,carneiroruggiero2,carneiroruggiero3,contrerasiturriagabook}.
In contrast, our result gives a sort of ``dual graph theorem'' for all minimizers, in the sense that, although the support of a  minimizer $\mu$ need not be contained in the graph of a section of the tangent bundle $TM$, the momenta $\partial L/\partial v$ on the support of $\mu$ do need to be contained in a section of the cotangent bundle $T^*M$. 

The literature is extensive for results giving sufficient conditions for the existence of minimizers, and these conditions typically come in the form of coercivity, boundedness, or super-linearity of the function $L$; see for example \cite{dacorogna,giusti}. In this direction there is also the line of research initiated by Morrey \cite{morrey1952quasi} related to quasi-convexity, a condition that was found to be necessary and sufficient for the weak sequential lower continuity of the action. %
In this paper, however, we assume the existence of the minimizers, and we do not examine the question of %
the continuity of the action.

We give examples of the application and sharpness of Theorem \ref{thm:mather} in Section \ref{sec:examples}. In particular, we show that, in the level of generality that we need to work in so as to obtain a full characterization, it is impossible to prove any regularity of the minimizers (Example \ref{ex:exactform}), and instead only the regularity of the momenta $df=\partial L/\partial v$ can be established (which is done in Section \ref{sec:lipschitzregularity}), and the regularity we prove is sharp (Examples \ref{ex:discontinuousdf}, \ref{ex:sharpregularity}, and \ref{ex:arnaud}). In particular, Bernard's result \cite{bernard2007existence} of existence of $C^{1,1}$ subsolutions of the Hamilton-Jacobi equation for Tonelli Lagrangians cannot be recovered in this setting. However we are able to recover Mather's original Lipschitz regularity result (Example \ref{ex:tonelli}). We also show how to use our theorem to prove some regularity of the distance function for non-strictly convex Finsler metrics (Example \ref{ex:finsler}).

A question that remains open is that of the regularity of the form $df$ at the boundary $\partial \mu$ of the minimizers.

\subsection{Optimal mass transport}

The optimal mass transport context can be formulated as follows \cite[\S7.2]{ambrosiogiglisavare}: given two probability measures $\nu_1$ and $\nu_2$ on a manifold $M$, the problem is to find a measure $\mu$ on $TM$ whose boundary (understood as the boundary of the current on $M$ induced by $\mu$) is  $\partial \mu=\nu_2-\nu_1$ and $\mu$ minimizes a Lagrangian action. The measure $\mu$ is understood to be encoding a bunch of curves that dictate how the mass of $\nu_1$ should be moved to where the mass of $\nu_2$ is, and the integrand is usually the arclength of those curves encoded by $\mu$, so that the overall interpretation is that one is finding the way to move $\nu_1$ into $\nu_2$ with the least possible effort. 

The Young Superposition Principle (Example \ref{ex:youngsuperposition}) is the statement that a measure $\mu$ with boundary $\partial\mu=\nu_2-\nu_1$ is a convex combination of generalized curves delineating trajectories that join $\nu_1$ and $\nu_2$, so that the description above is meaningful.

There are a few formulations of the optimal mass transport problem, according to the properties of the Lagrangian function $L$ involved in the action. 

For the case in which the Lagrangian $L$ is a $C^2$ function throughout the tangent bundle,
our results in Theorem \ref{thm:mather} and Corollary \ref{cor:mathertimedependent} characterize the Lagrangians for which this problem can be considered meaningfully, draw a relation with the corresponding Hamilton--Jacobi equation, and give \emph{a priori} $C^{1,1}$-regularity on the projected support of $\mu$ for the subsolutions of that equation that govern the properties of $\mu$; cf. \cite[Theorem 6.2.7]{ambrosiogiglisavare}. Our results also generalize the Lipschitz regularity result obtained by Bernard--Buffoni \cite{bernardbuffoni}, where the authors rely on assumptions of convexity, super-linearity, and completeness of the geodesic flow, which are not necessary in our version.

Another formulation of the optimal transport problem involves Lagrangian densities $L$ that fail to be $C^2$ at the origin, but instead are positively homogeneous of degree 1, namely, $L(x,\lambda v)=\lambda L(x,v)$ for all $\lambda>0$, $(x,v)\in TM$, such as is the case for Finsler metrics. In this case, Lipschitz regularity results were also obtained by Bernard--Buffoni \cite{bernardbuffonimonge} and other authors mentioned in their work, and our results in Example \ref{ex:finslertransport} generalize those, as we manage to avoid the convexity assumption. Also in that example, in the case in which $L$ is the distance associated to a Riemannian metric, we recover the result \cite{ambrosiogiglisavare} that the transport is done along a gradient flow.

\subsection{Optimal control}

We also apply the same line of reasoning developed for Theorem \ref{thm:mather} to analyze an optimal control problem in Section \ref{sec:ocontrol}, and we are able to give a coarse characterization of minimizable integrands in Theorem \ref{thm:ocontrolcharact} akin to the one developed in Section \ref{sec:coarsecharact}, as well as a result, Theorem \ref{thm:matherocontrol}, that gives sufficient conditions to obtain Lipschitz regularity of the momenta. We explain in Remark \ref{rmk:HJB} the close connection with the Hamilton--Jacobi--Bellman equation satisfied by the value function and the Maximum Principle.

\subsection{Other remarks}

The author's view is that this paper is about the importance of the holonomy constraint $\partial\mu=c$ (see Definition \ref{defn:currents}), from whose exploitation stem most of the results that we obtain here. A recurrent motif is that whenever we minimize within a set of measures that vanish for a certain class of functions (exact forms in our case), the minimizable functionals correspond to functions that are nonnegative up to a addition of function in that class. What is interesting then is how the seemingly innocent assumption of minimizability in a set that satisfies $\partial\mu=c$ implies already some regularity as well as familiar concepts like energy conservation, the existence of calibrations, the maximum principle in optimal control, and the ubiquity of the Hamilton--Jacobi equation, among others.

The results of this paper have been applied to the characterization of deformations of closed measures in \cite{mydeformations}.

\section{Characterization of minimizable action functionals}
\label{sec:characterization}
\subsection{Coarse characterization}\label{sec:coarsecharact}
\subsubsection{Time-independent setting}

Let $d$ be a positive integer, and let $M$ be a $d$-dimensional, second-countable, smooth manifold. Denote by $TM$ its tangent bundle, with fibers $T_xM\cong{\mathbb{R}}^d$ at each point $x$ in $M$.  

 On a manifold $X$,
 let $C^\infty(X)$ denote the set of smooth functions on $X$ with the topology induced by the seminorms $|\cdot|_{K,k}$ associated to each compact subset $K$ of $X$ and to each positive integer $k$, and given by 
 \begin{equation}\label{eq:seminorms}
  |f|_{K,k}=\sum_{|I|\leq k}\sup_{x\in K}\left|\partial^If(x)\right|. 
 \end{equation}

\begin{definition}\label{defn:currents}
 A \emph{(compactly-supported) normal 0-current} is a continuous real-valued functional $C^\infty(M)\to {\mathbb{R}}$ given by a signed, Radon, compactly-supported measure on $M$. %
 Given a compactly-supported, Radon measure $\mu$ on $TM$, the \emph{boundary} \[\partial \mu\colon C^1(M)\to{\mathbb{R}}\] of $\mu$ is defined by
 \[\langle\partial\mu,f\rangle=\int_{TM}df\,d\mu,\quad f\in C^1(M).\]
 A measure $\mu$ is \emph{closed} if $\partial\mu=0$. It can be checked that, if $M$ is connected, a normal 0-current $c$ is a boundary if, and only if, $\langle c,1\rangle=0$.

 For a fixed normal 0-current $c$, let $\mathscr H(c)$ be the set of compactly-supported, positive, Radon measures $\mu$ with $\partial \mu=c$; in the special case $c=0$, we additionally require the elements of $\mathscr H(0)$ to be probability measures.
\end{definition}

Let $E$ be a complete, sequential, locally-convex topological vector space of Borel measurable functions on $TM$ that contains $C^\infty(TM)$ as a subset. We will assume that $\mathscr H(c)\subset E^*$ and that the topology of $C^\infty(TM)$ described above is finer than the one this space inherits from $E$, so that every open set in the inherited topology is an open set in the topology described above. This assumption implies that every continuous linear functional $\vartheta\in E^*$ defines a compactly-supported distribution when restricted to $C^\infty(TM)$.

To give some examples, $E$ could be %
a space of $C^\ell$ functions on $TM$ with $\ell\in[0,\infty]$ and the topology induced by the seminorms \eqref{eq:seminorms} for $k\leq\ell$.
For the verification of their adequacy, it may be useful to recall that it is enough for a topological vector space to be normed, metric, or 
first-countable in order for it to be sequential.

\begin{theorem}\label{thm:characterization}
 Let $c$ be a normal 0-current on $M$ that is a boundary. 

 If $L$ is an element of $E$ such that the action functional $\nu\mapsto\int_{TM} L\,d\nu$ reaches its minimum within $\mathscr H(c)$ at some point $\mu$, then 
 there exist functions $f_1,f_2,\dots$ in $C^\infty(M)$, nonnegative functions $g_1,g_2,\dots$ in $E$ such that 
 \[\lim_{i\to+\infty}\int_{TM} g_i\,d\mu=0,\]
and
 \[L=\begin{cases}
      \int_{TM} L\,d\mu+\lim_{i\to+\infty}(df_i+g_i),&c=0,\\
      \lim_{i\to+\infty}(df_i+g_i),&c\neq0,
     \end{cases}
\]
 where the limits are taken in $E$.  In particular, if $c\neq 0$,
 \[\int L\,d\mu=\lim_i\langle c,f_i\rangle.\] 
 Conversely, if $L$ has this structure, its action reaches its minimum within $\mathscr H(c)$ at $\mu$.
\end{theorem}

We immediately have the following consequence of the fact that $\int g_id\mu\to0$.

\begin{corollary}[Energy conservation]\label{coro:energyconservation}
 Let $L$ be an element of $E$ and %
 assume that the action functional of $L$ reaches its minimum within $\mathscr H(c)$ at a minimizer $\mu$. Define the \emph{Hamiltonian} associated to $L$ to be the function on $T^*M$ given by
 \[H(x,u)=\sup_{v\in T_xM}u_x(v)-L(x,v)\quad \textrm{for}\quad x\in M, u\in T^*_xM.\]
 Then the value of $H$ is constant throughout $\mu$-almost all the support of $\mu$ and equals 
 \[
  H(x,u)=\begin{cases}-\int L\,d\mu,&c=0,\\
            0,&c\neq 0,
          \end{cases}
          \quad 
          u=\frac{\partial L}{\partial v}(x,v),\;(x,v)\in\operatorname{supp}\mu.
\]
\end{corollary}

\begin{remark}[Terminology]\label{rmk:energyconservationterminology}
 The statement of Corollary \ref{coro:energyconservation} is known as an \emph{Energy Conservation Principle} because of historical reasons, which we now sketch. The quantity $H$ is known as the \emph{Hamiltonian energy}. If we add a time coordinate (of which $L$ remains independent) and situate ourselves in the time dependent setting (see Section \ref{sec:timedependent}), then the minimizing measure $\mu$ satisfies \cite[Lemma 16]{patrick}, in the sense of distributions, a continuity equation of the form
 \[\partial_t \mu+\operatorname{div}(V\mu)=0,\]
 where $V\colon TM\times {\mathbb{R}}\to TM$ is the average velocity on each fiber $T_xM\times\{t\}$, which allows one to interpret $\mu$ as flowing along the vector field $V$. A particular consequence of the corollary is that the quantity $H$ remains constant for almost all times $t$, and this gives meaning to the word ``conservation'' in the present context. (Theorem \ref{thm:mather} will show that, with slightly stronger assumptions on $L$, $H$ in fact is constant for all $t$.)
\end{remark}

For the proof of Theorem \ref{thm:characterization}, we will need the following lemma.

\begin{lemma}\label{lem:dual}
 In the setting of Theorem \ref{thm:characterization}, let
 \begin{align*}
   &Q=\{\ell\colon E^*\to\mathbb{R} \mid \textrm{$\ell$ is affine and continuous, its linear part is induced by}\\
   &\qquad\quad\textrm{evaluation at an element of $E$, and $\ell(\mu)\geq 0$ for all $\mu\in \mathscr H(c)$}\},\\
   &R=\{\ell\colon E^*\to\mathbb R\mid \ell(\xi)=\xi(df+g)-\langle c,f\rangle,\;  f\in C^\infty(M),\;g\in E,\;g\geq 0 \}.
 \end{align*} 
 Then we have $Q=\overline R$ in $E\oplus \mathbb R\subseteq E^{**}\oplus\mathbb R$.
\end{lemma}
For the proof of Lemma \ref{lem:dual}, we recall
\begin{definition}
 Let $V$ be a locally-convex topological vector space with topological dual $V^*$. For a set $A\subseteq V$, the \emph{dual} $A'\subseteq V^*$ is the set consisting of all affine functionals $\theta\in V^*\oplus\mathbb R$ such that $\theta(v)\geq 0$ for all $v\in A$.
\end{definition}
\begin{remark}\label{rmk:duality}
 It is an easy consequence of the Hahn--Banach Separation Theorem that 
 if $X$ and $Y$ are two convex cones in a locally-convex topological vector space $V$ and $X'=Y'$, then their closures coincide, $\overline X=\overline Y$.
\end{remark}

\begin{proof}[Proof of Lemma \ref{lem:dual}]
 We will show that $Q'=R'$ in $E^*\oplus\mathbb R$; in fact, these dual sets coincide with $\mathbb R_{\geq 0}\mathscr H(0)$ when $c=0$ and with $\mathscr H(c)$ when $c\neq 0$.
 
 To see why, note first that the set of linear functionals $\ell_g(\xi)=\xi(g)$ induced by nonnegative elements $g\in C^\infty (TM)$, $g\geq0$, is a subset of both $Q$ and $R$. So if $\xi\in Q'$ or $\xi\in R'$, necessarily $\ell_g(\xi)\geq 0$ for all $g$, and by \cite[\S6.22]{liebloss} $\xi$ can be represented as integration over a compactly supported, nonnegative, Radon measure.
 
 Also, if $f\in C^\infty(M)$, then the affine functional 
 \[\ell_{f}(\xi)=\xi(df)-\langle c,f\rangle\]
 belongs to both $Q$ and $R$, so it is nonnegative throughout $Q'$ and $R'$. But this is also true if we replace $f$ by $-f$, so we also have $0\leq\ell_{-f}(\xi)=-\ell_f(\xi)$ for all $\xi\in Q'\cup R'$. We conclude that 
 \[0=\ell_f(\xi)=\xi(df)-\langle c,f\rangle=\partial \xi(f)-\langle c,f\rangle,\quad\xi\in Q'\cup R'.\]
 This shows that $\xi$ is contained in $\mathbb R_{\geq 0}\mathscr H(0)$ if $c=0$ or in $\mathscr H(c)$ if $c\neq 0$.
\end{proof}

\begin{proof}[Proof of Theorem \ref{thm:characterization}]
 Take a function $L\in E$ that satisfies
 \[\int L\,d\mu\leq\int L\,d\nu\quad\textrm{for all}\quad\nu\in\mathscr H(c),\]
 and consider the affine functional on $E^*$ given by
 \[\ell(\xi)=\langle \xi, L\rangle -\int L\,d\mu.
\] Then $\ell$ satisfies
 \[0=\ell(\mu)\leq \ell(\nu)\quad\textrm{for all}\quad\nu\in\mathscr H(c).\]
 Thus $\ell$ belongs to the set $Q$ in the statement of Lemma \ref{lem:dual}. It follows that $\ell$ also belongs to $\overline R$. Since  $E$ is sequential, the topological closure equals the sequential closure, so there exists a sequence of functionals $\ell_i\in R$ of the form
 \[\ell_i(\xi)=\xi(df_i+g_i) -\langle c,f_i\rangle,\quad i=1,2,\dots,\] 
 converging to $\ell_i\to\ell$ and such that  $f_i\in C^\infty(M)$, $g_i\in E$, $g_i\geq 0$. 
 
 Comparing the linear and constant parts of $\ell$ and $\ell_i$, we conclude that $\int L\,d\mu=\lim_i\langle c,f_i\rangle$ if $c\neq 0$.

 We also have that 
 \[0=\ell(\mu)=\lim_i\ell_i(\mu)=\lim_i\int g_i\,d\mu\]
 since $\int df_i\,d\mu=\langle c,f_i\rangle$ by definition.
\end{proof}

\subsubsection{Time-dependent setting}\label{sec:timedependent}
\begin{definition}\label{defn:timedependent}
 When talking about the \emph{time-dependent setting}, we will refer to the situation in which $M$ is of the form $N\times P$, where $N$ is a smooth $(d-1)$-dimensional manifold and $P$ is a connected, 1-dimensional, smooth manifold that plays the role of time. 
 We fix a parameterization $t\colon I\subseteq{\mathbb{R}}\to P$, and we use it to distinguish, at each point $p\in P$, the vector $\mathbf 1\in T_pP$ such that $dt_p\mathbf 1=1$.
 We will denote by $\mathscr H_1(c)\subset\mathscr H(c)$ the set of elements of $\mathscr H(c)$ %
 that are supported within the set $TN\times P\times\{\mathbf 1\}\subset TM$, which we will identify with $TN\times P$. 
\end{definition}
\begin{remark}
 We interpret measures $\mu$ in $\mathscr H_1(c)$ as advancing in the time direction with ``velocity'' $\mathbf 1\in T_pP\cong {\mathbb{R}}$, $p\in P$, which roughly means that time itself always moves a the same speed.
\end{remark}

Let $E$ be a complete, sequential, locally-convex topological vector space of Borel measurable functions on $TN\times P$ that contains $C^\infty(TN\times P)$, such that $\mathscr H_1(c)\subset E^*$, and such that the topology inherited by $C^\infty(TN\times P)$ from $E$ is finer than the topology induced by the seminorms \eqref{eq:seminorms}. For example, $E$ could be %
$C^\ell(TN\times P)$ for $\ell\in[0,\infty]$ with the topology induced by the seminorms \eqref{eq:seminorms} for $k\leq \ell$.

It is straightforward to adapt the reasoning that gives Theorem \ref{thm:characterization} in order to obtain
\begin{corollary}\label{coro:timedepcoarse}
 Let $N\times P$ be a manifold that is the product of $C^\infty$ manifolds $N$ and $P$, with $P$ playing the role of time, $\dim P=1$, so that we are in the time-dependent setting. Let $c$ be a normal 0-current on $N\times P$ such that $\mathscr H_1(c)$ is not empty.
 Assume that $L$ is an element of $E$ such that $\nu\mapsto\int_{TN\times P} L\,d\nu$ reaches its minimum within $\mathscr H_1(c)$ at some point $\mu$. Then there exist functions $f_1,f_2,\dots$ in $C^\infty(N\times P)$, nonnegative functions $g_1,g_2,\dots$ in $E$ such that 
 \[
 \lim_{i\to+\infty}\int_{TM} g_i\,d\mu=0,\]
 and
 \[L=\lim_{i\to+\infty}(d_{TN}f_i+\frac{\partial f_i}{\partial t}+g_i),\]
 where the limit is taken in $E$, %
 and $d_{TN}f_i\colon TN\times P\to{\mathbb{R}}$ indicates the exterior derivative of the restriction of $f_i$ to $TN\times\{t_0\}$ for the corresponding point $t_0\in P$.
 
 Conversely, if $L$ has this structure, its action reaches its minimum within $\mathscr H_1(c)$ at $\mu$.
\end{corollary}

This follows from Theorem \ref{thm:characterization} taking $M=N\times P$.
Note that $\partial f_i/\partial t$ is constant on each fiber $T_xN$.

Although a version of the energy conservation result, Corollary \ref{coro:energyconservation}, holds in the time-dependent setting, it is not very transparent because it involves the less-tangible $\partial f_i/\partial t$.

\subsection{Lipschitz regularity}\label{sec:lipschitzregularity}
\subsubsection{Time-independent setting}

\begin{definition}
 Let $x$ be a point in an open set $U\subseteq M$.
 The form $\theta\in T_x^*M$ is a \emph{Clarke differential of $f\colon U\to {\mathbb{R}}$ at $x$} if it is in the convex hull of the accumulation points of the values of $df_y$ as $x\to y$. A section $\alpha\colon U\to T^*U$ is a \emph{Clarke differential of $f\colon U\to{\mathbb{R}}$} if it is a Clarke differential at every point of $U$.
\end{definition}

Denote by $\pi\colon TM\to M$ the fiberwise projection.

\begin{theorem}\label{thm:mather}
 Let $c$ be a normal 0-current on $M$ that is a boundary, so that the space $\mathscr H(c)$ from Definition \ref{defn:currents} is nonempty. 
 Let $L$ be an element of $C^2(TM)$ such that the action functional $\nu\mapsto\int_{TM} L\,d\nu$ reaches its minimum within $\mathscr H(c)$ at some point $\mu$. Let $U$ be an open subset of $M$ with compact closure and such that $\pi(\operatorname{supp}\mu)\subseteq U$. Then there exist a Lipschitz function $f\colon U\to {\mathbb{R}}$, a nonnegative function $g\colon TU\to {\mathbb{R}}_{\geq 0}$, and a bounded (possibly discontinuous) section $\alpha\colon U\to T^*U$ such that:
 \begin{enumerate}
  \item \label{it:alphaclarke} $\alpha$ is a Clarke differential of $f$ (coinciding with $\alpha=df$ wherever $f$ is differentiable);
  \item \label{it:sumdecomp} throughout $TU$, we have
   \[L=\begin{cases}\int_{TM} L\,d\mu+\alpha+g,&c=0,\\\alpha+g&c\neq 0,\end{cases}\]
   which in particular means that 
   \[L=\begin{cases}\int_{TM} L\,d\mu+df+g,&c=0,\\ df+g&c\neq 0,\end{cases}\] 
   wherever $f$ is differentiable;
  \item \label{it:differentiability} $f$ is differentiable on $\pi(\operatorname{supp}\mu)\setminus\operatorname{supp} c$ and, for $(x_0,v_0)\in\operatorname{supp}\mu$ with $x_0\notin\operatorname{supp} c$, we have
   \[df_{x_0}=\frac{\partial L}{\partial v}(x_0,v_0);\]
  \item \label{it:lipschitzity} for every open neighborhood $Y\subset U$ of $\operatorname{supp} c$, there is a constant $C_Y>0$ such that the map $x\mapsto df_x$ is $C_Y$-Lipschitz throughout $\pi(\operatorname{supp}\mu)\cap U\setminus Y$;
  \item \label{it:gvanishes} $g\equiv 0$ throughout $\operatorname{supp}\mu$.
 \end{enumerate}
 Conversely, if $L$ has this structure, its action reaches its minimum within $\mathscr H(c)$ at $\mu$.
\end{theorem}
\begin{remark}
 In items \ref{it:differentiability} and \ref{it:lipschitzity} we work away from the support of $c$ for simplicity. Our proof shows however that the Lipschitzity should hold on all points $x\in\pi(\operatorname{supp}\mu)$ for which there exist absolutely continuous $\gamma\colon[-\varepsilon,\varepsilon]\to M$ with $\varepsilon>0$, $\gamma(0)=x$, and $\gamma'(t)\subseteq\operatorname{conv}(\operatorname{supp}\mu\cap{T_{\gamma(t)}M})$ for almost every $t\in[-\varepsilon,\varepsilon]$.
\end{remark}

\begin{proof}
 By replacing $L$ with $L-\int L\,d\mu$, we may assume that $\int L\,d\mu=0$.

 By Theorem \ref{thm:characterization}, there are functions $f_1,f_2,\dots$ in $ C^\infty(M)$ and $g_1,g_2,\dots$ in $ C^2(TM)$ such that
 \[L=\lim_{i\to+\infty}g_i+df_i\]
 in $C^2(TM)$ with the topology induced by the seminorms \eqref{eq:seminorms} with $k\leq 2$, $g_i\geq 0$, $\lim_{i\to\infty}\langle c,f_i\rangle=0$, 
 and $\lim_{i\to\infty}\int g_i\,d\mu=0$. We may additionally assume that $f_1,f_2,\dots$ are uniformly bounded on $\overline U$.
 
 Since both $f_i$ and $df_i$ are uniformly bounded on the compact set $\overline U$, by an application of Arzel\`a-Ascoli, perhaps passing to a subsequence, we may assume that the sequence $\{f_i\}_i$ converges on $C^0(\overline U)$ to a Lipschitz function $f\colon U\to M$. By Rademacher's theorem, $f$ is differentiable almost everywhere on $U$. Since all $df_i$ satisfy $df_i\leq L_i$ and $L_i\to L$, it can be checked that $df\leq L$, wherever $df$ is defined. 
 
 We let $\alpha$ be a Clarke differential of $f$ on $U$; since $df\leq L$ wherever it is defined, by continuity of $L$ we know that it is possible to choose $\alpha$ so that $\alpha\leq L$. We set $g=L-\alpha\geq 0$ on $TU$. With these definitions, we have that items \ref{it:alphaclarke} and \ref{it:sumdecomp} hold.
 
 Similarly, since 
 \[\int g\,d\mu=\lim_i\int g_i\,d\mu=0\] and $g\geq0$, it follows that $g=0$ $\mu$-almost everywhere. It will follow from the continuity of $df$ on $\operatorname{supp}\mu$ (a consequence of item \ref{it:lipschitzity}) that $g$ actually vanishes throughout $\operatorname{supp}\mu$, as stated in item \ref{it:gvanishes}. 
 
 It remains to show that items \ref{it:differentiability} and \ref{it:lipschitzity} are true. We will prove these items for each open subset $V$ of $U$ that does not intersect $\operatorname{supp} c$ and is diffeomorphic to the open ball in ${\mathbb{R}}^d$, and the result will follow from the compactness of $\operatorname{supp}\mu$. We pass through the chart of $V$ to the unit ball in ${\mathbb{R}}^d$, but for simplicity we keep all notations the same.
 
 Let $K$ be a subset of $TV$ of the form $\{(x,v)\in TV:|v|<C_0\}$ for some $C_0\gg0$ such that $\operatorname{supp}\mu\cap TV\subset K$. Let $\phi\colon TV\to{\mathbb{R}}_{\geq 0}$ be a smooth function that vanishes on $K$ and grows larger than a positive multiple of $|v|^2$ outside a neighborhood of $K$. Note that replacing $L$ with $L+\phi$ changes neither $f$, nor $df$, nor the statements of items \ref{it:differentiability} and \ref{it:lipschitzity}; thus, we may assume, without loss of generality, that the Lagrangian $L$ grows super-quadratically in $v$.

 The convexification $\tilde L$ of $L$ is defined by 
 \begin{equation}\label{eq:convexification}
  \tilde L(x,v)=\sup\{r+\theta(v):r\in {\mathbb{R}}, \theta\in T^*_xV,r+\theta\leq L|_{T_xV}\}, \quad (x,v)\in TV.
 \end{equation}
 Since $L$ is super-quadratic in $v$, $\tilde L$ is super-quadratic too. It is locally Lipschitz throughout $V$. Also, it follows from \cite[Theorem 4.2]{griewankrabier} that $\tilde L|_{T_xV}$ is in $C^{1,1}_{\operatorname{loc}}(T_xV)$ for each $x\in V$, that is, it is differentiable and its derivative is locally Lipschitz. We show in Lemma \ref{lem:Lvlipschitz} that $(x,v)\mapsto\frac{\partial \tilde L}{\partial v}(x,v)$ is locally Lipschitz.
 The Lagrangian $\tilde L$ can also be decomposed as
 \[\tilde L=\alpha+\tilde g\]
 for some function $0\leq \tilde g\leq g$ that is convex and $C^{1,1}_{\operatorname{loc}}$ in the fibers of $TV$. 
 
 The rest of the proof is inspired in \cite[Section 4.11]{fathibook}. It follows from Lemma \ref{lem:smirnov} that for $\pi_*\mu$-almost every $x$ in $\pi(\operatorname{supp}\mu)\cap V$ there is a curve $\gamma_x$ as in the hypothesis of Lemma \ref{lem:lipcondholds}; let us call $A$ the set of such points $x$. This means that the hypothesis for the criterion for $f$ having a Lipschitz derivative, Lemma \ref{lem:lipschitzcriterion}, hold for all points in $A\cap V$, so that $f$ is differentiable throughout $A$ and its derivative $x\mapsto df_x$ is locally Lipschitz in the set $A_{\frac13}=A\cap\frac13 \overline V$. 
 
While $x\mapsto\dot\gamma_x(0)$ may not be continuous, the family of linear functionals
 \[
  A_{\frac13}\ni x\mapsto \frac{\partial \tilde L}{\partial v}(x,\dot\gamma_x(0))\in T^*_xV
 \]
is Lipschitz, so it can be extended in a unique way to the closure $\overline A_{\frac13}$. Additionally, since $\tilde L$ and $\partial \tilde L/\partial v$ are locally Lipschitz, it follows that at each point $x$ in $\overline A_{\frac13}$ there is a vector $v_x$ such that the Lipschitz extension can be written as
\[
  \overline A_{\frac13}\ni x\mapsto \frac{\partial \tilde L}{\partial v}(x,v_x)\in T^*_x(\tfrac13V),
\]
with $v_x=\dot\gamma_x(0)$ for $x\in A_{\frac13}$.
(It may be that the map $x\mapsto v_x$ is discontinuous however.)

The functions $f$ and $x\mapsto \partial\tilde L/\partial v(x,v_x)$ being locally Lipschitz (by Lemma \ref{lem:Lvlipschitz}), the condition 
 \begin{equation*}
  |f(y)-f(x)-\frac{\partial \tilde L}{\partial v}(x,v_x)|\leq K\|y-x\|^2
 \end{equation*}
 must hold throughout the closure $\overline A_{\frac13}$ (replacing $\dot\gamma(0)$ with $v_x$). By another application of Lemma \ref{lem:lipschitzcriterion}, it follows that $f$ is actually differentiable throughout $\overline {A\cap\frac19V}\supset \pi(\operatorname{supp} \mu)\cap \frac19V$ and that the map $x\mapsto df_x$ is Lipschitz there, as we wanted to prove.
\end{proof}

\begin{lemma}\label{lem:Lvlipschitz}
 Let $A$ and $B$ be two open, convex subsets of ${\mathbb{R}}^m$ and ${\mathbb{R}}^n$, respectively, for some $m,n>0$.
 Assume that the function $c\colon A\times B\to{\mathbb{R}}$ is Lipschitz, and for each $u\in A$ the map $v\mapsto c(u,v)$ is convex and $C^{1,1}_{\operatorname{loc}}$. Then $\partial c/\partial v$ is locally Lipschitz.
\end{lemma}
\begin{proof}
  Since $c$ is known to be $C^{1,1}_{\operatorname{loc}}$ in the $B$-direction, we need only show that, for each fixed vector $v_0$, the mapping $u\mapsto\partial c/\partial v(u,v_0)$ is locally Lipschitz.
  
  For $v_1,v_2,\dots,v_{n+1}\in B$ be points in general position, and let 
  \[f_{v_1,\dots,v_{n+1}}\colon A\times B\to{\mathbb{R}}\] 
  be the function, linear in the second coordinate $v\in B$, such that 
  \[f_{v_1,\dots,v_{n+1}}(u,v_i)=c(u,v_i).\] 
  Note that, for fixed $u\in A$ and for $v$ within the simplex generated by $v_1,\dots, v_{n+1}$, we have $f_{v_1,\dots,v_{n+1}}(u,v)\geq c(u,v)$ by the convexity of $v\mapsto c(u,v)$.
  For each neighborhood $W\subseteq B$ of $v_0$, let $\Gamma_W(u)$ be the set of all vectors $\frac{\partial f_{v_1,\dots,v_{n+1}}}{\partial v}(u,v_0)$ for $v_1,\dots,v_{n+1}$ such that $v_0$ is in the simplex they generate. Then it is easy to see (for example by applying the mean value theorem in the directions of the standard basis of ${\mathbb{R}}^n$ and using the fact that $v\mapsto c(u,v)$ is $C^{1,1}_{\operatorname{loc}}$) that 
  \[\frac{\partial c}{\partial v}(u,v_0)=\bigcap_W\Gamma_W(u),\]
  where the intersection is taken over all neighborhoods $W\subseteq B$ of $v_0$.
  Since $c$ is Lipschitz, the maps $u\mapsto f_{v_1,\dots,v_{n+1}}(u,v_i)$ are Lipschitz too, and so also the maps $u\mapsto \frac{\partial f_{v_1,\dots,v_{n+1}}}{\partial v}(u,v_0)$ must be Lipschitz, uniformly so for $v_1,\dots,v_{n+1}$ in small neighborhoods $W$ around $v_0$. Hence the contents of each $\Gamma_W(u)$ varies in a Lipschitz way with $u$, as does their intersection. This proves the lemma.
\end{proof}

\begin{lemma}\label{lem:Lsuperdiff}
 Fix an open subset $D$ of $TV$ with compact closure $\overline D$.
 Then there is a constant $C_D>0$ such that, for every $(x_0,v_0)\in \overline D$ and for every smooth section $\sigma\colon \overline V\to D$ with $\sigma(x_0)=v_0$, there is $\eta\in T_x^*V$ such that
 \[\tilde L(x,\sigma(x))-\tilde L(x_0,v_0)-\eta(x-x_0)\leq C_D|x-x_0|^2\]
 for all $x$ in a neighborhood of $x_0$.
 In particular, $\tilde L\circ \sigma$ is everywhere upper Dini differentiable in $V$. 
\end{lemma}
\begin{proof}
 Let $C_D$ be an upper bound for the norm of the Hessian of $L$ in $\overline D$.

 For each $(d+1)$-tuple of vectors $v_1,v_2,\dots,v_{d+1}\in T_{x_0}V\cong {\mathbb{R}}^d$ and for each $x\in V$ such that $\sigma(x)$ is in the interior of their convex hull, let 
 \[\lambda_1(v_1,\dots,v_{d+1};x),\dots,\lambda_{d+1}(v_1,\dots,v_{d+1};x)\geq0\]
 be such that $\sum_i\lambda_i(v_1,\dots,v_{d+1};x)=1$ and $\sum_i\lambda_i(v_1,\dots,v_{d+1};x)v_i=\sigma(x)$. Also let
 \[\phi(v_1,\dots,v_{d+1};x)=\sum_{i=1}^{d+1}\lambda_i(v_1,\dots,v_{d+1};x)L(x,v_i).\]
 Then if $v_0$ is in the interior of $\operatorname{conv}(v_1,\dots,v_{d+1})$, $x\mapsto \phi(v_1,\dots,v_{d+1};x)$ is defined and $C^2$ in a neighborhood of $x_0$ and 
 \[\tilde L(x,\sigma(x))=\inf\{\phi(v_1,\dots,v_{d+1};x):v_0\in\operatorname{conv}(v_1,\dots,v_{d+1})\}.\]
 Since the second derivatives of all the functions on $\phi(v_1,\dots,v_{d+1};\cdot)$ are bounded by $C_D$, the statement of the lemma follows.
\end{proof}

\begin{lemma}\label{lem:smirnov}
 For $\pi_*\mu$-almost every $x\in\pi(\operatorname{supp}\mu)\cap V$ there is some $t_0>0$ such that for almost every $0<t\leq t_0$, there is an absolutely continuous curve $\gamma\colon[-t,t]\to \pi(\operatorname{supp}\mu)\cap V$ such that
 \begin{enumerate}
 \item  $\gamma(0)=x$, 
 \item for almost every $s\in[-t,t]$ the velocity $\dot\gamma(s)$ is defined and contained in the convexification of $\operatorname{supp}\mu\cap T_{\gamma(s)}V$,  
 \item for all $-t\leq a\leq b\leq t$ we have $\int_{a}^{b} \tilde L(d\gamma(s))ds=f(\gamma(b))-f(\gamma(a))$, and 
 \item $\gamma$ minimizes the action of the convexified Lagrangian $\tilde L$ among all absolutely continuous curves with endpoints $\gamma(a)$ and $\gamma(b)$.
 \end{enumerate}
\end{lemma}
\begin{proof}
 This follows immediately from the decomposition result of \cite{smirnov1993decomposition}; see also the exposition in \cite[\S 3]{bangert1999minimal}.
\end{proof}

\begin{corollary}\label{coro:laxoleinik}
 For  $\pi_*\mu$-almost every $x\in \pi(\operatorname{supp}\mu)\cap V$ there is $t_0>0$ such that for all $0<t\leq t_0$, the function $f$ satisfies
 \begin{equation}\label{eq:laxoleinikpast}
  f(x)=\inf_\gamma f\circ\gamma(-t)+\int_{-t}^0 \tilde L(d\gamma)\,ds
 \end{equation}
 where the infimum is taken over all absolutely-continuous curves $\gamma\colon[-t,0]\to V$ with $x=\gamma(t)$.
 
 For  $\pi_*\mu$-almost every $x\in \pi(\operatorname{supp}\mu)\cap V$ there is $t_0>0$ such that for all $0<t\leq t_0$, the function $f$ satisfies
 \begin{equation}\label{eq:laxoleinikfuture}
  f(x)=\sup_\gamma f\circ\gamma(t)-\int_0^t \tilde L(d\gamma)\,ds
 \end{equation}
 where the supremum is taken over all absolutely-continuous curves $\gamma\colon[0,t]\to V$ with $x=\gamma(0)$.
 
 Moreover, the infimum \eqref{eq:laxoleinikpast} and the supremum \eqref{eq:laxoleinikfuture} are both realized by ab\-so\-lute\-ly-continuous curves whose images are contained in $\pi(\operatorname{supp}\mu)$.
\end{corollary}

\begin{lemma}\label{lem:existenceLx}
 Let $\gamma\colon[-t,t]\to V$ be an absolutely continuous curve passing through $\gamma(0)=x\in\pi(\operatorname{supp}\mu)$ and such that $\gamma$ minimizes the action of $\tilde L$ among all absolutely continuous curves with the same endpoints $\gamma(-t)$ and $\gamma(t)$. Then for almost every $s\in[-t,t]$, $d\gamma(s)$ is defined in the sense that $s$ is a Lebesgue point of $d\gamma\in L^1([-t,t])$, and we have that $\tilde L$ is differentiable at $d\gamma(s)$.
 
 Also, for all smooth $h\colon [-t,t]\to {\mathbb{R}}^d$ with $h(-t)=0=h(t)$ %
 we have
 \begin{equation}\label{eq:nullonvanishing}
  \int_{-t}^t \tilde L_x(d\gamma(s))h(s)+\tilde L_v(d\gamma(s))\dot h(s)\,ds=0.
 \end{equation}
\end{lemma}
\begin{proof}
 The fact that $d\gamma$ is defined almost everywhere on $[-t,t]$ follows from the Lebesgue Differentiation Theorem.
 We know that $\tilde L$ is $C^{1,1}_{\operatorname{loc}}$ in the fibers, so in order to prove the first statement we need only show that the derivative exists in the horizontal direction, that is, we need to show the existence of $\tilde L_x(d\gamma(s)$ for almost every $s\in[-t,t]$.
 
 Here we use the ``little-$o$ notation'', so that $o(1)$ and $o(\delta)$ stand for functions such that $o(1)\to0$ and $o(\delta)/\delta\to0$ as $\delta\to0$, respectively.
 
 Recall that, by virtue of Lemma \ref{lem:Lsuperdiff}, we know that for each point $(x,v)$ in $T{\mathbb{R}}^d$ there is a linear form $\eta_{(x,v)}$ such that, if $|q|=1$,
 \begin{equation}\label{eq:dominanceoflinearform}
  \eta_{x,v}(q)+o(\delta)\geq \frac{\tilde L(x+\delta q,v)-\tilde L(x,v)}{\delta},  
 \end{equation}
 as $\delta\searrow0$.
 
 Let $h\colon[-t,t]\to {\mathbb{R}}^d$ be a smooth map such that $h(-t)=0=h(t)$. Then we have, by virtue of \eqref{eq:dominanceoflinearform} and of the minimality of $\gamma$, and by a degree 1 Taylor expansion in the fiberwise direction, 
 \begin{align*}
  \int_{-t}^t &\eta_{d\gamma(s)} (h)ds=\frac1\delta\int_{-t}^t \eta_{d\gamma(s)} (\delta h)ds\\
  &\geq \frac1\delta\int_{-t}^t\tilde  L(\gamma+\delta h,\dot\gamma)-\tilde  L(\gamma,\dot\gamma)\,ds +o(1)\\
 \end{align*}
 \begin{align*} 
  &=\frac1\delta\int_{-t}^t\tilde  L(\gamma+\delta h,\dot\gamma+\delta \dot h)-\delta\tilde  L_v(\gamma+\delta h,\dot\gamma)\dot h+o(\delta)-\tilde L(\gamma,\dot\gamma)\,ds\\
  &\geq -\int_{-t}^t \tilde L_v(\gamma+\delta h,\dot \gamma)\dot h \,ds+o(1).
 \end{align*}
 As $\delta\searrow0$, this sandwiches the expression whose limit would correspond to the integral of the derivative in the horizontal direction (whose existence we want to prove),
 namely,
 \[\tilde L_x(\gamma,\dot\gamma)h=\lim_{\delta\searrow0} \frac1\delta\left(\tilde L(\gamma+\delta h,\dot\gamma)-\tilde  L(\gamma,\dot\gamma)\right),\]
 between two expressions that are linear in $h$. These expressions must hence coincide, thus implying the existence of $\tilde L_x$ for almost every $s\in[-t,t]$. Equation \eqref{eq:nullonvanishing} also follows immediately from this argument. 
\end{proof}

\begin{lemma}[{Criterion for a Lipschitz derivative  \cite[Proposition 4.11.3]{fathibook}}]\label{lem:lipschitzcriterion}
 Let $B$ be the open unit ball in the normed space $E$.
 Fix a map $h\colon B\to {\mathbb{R}}$. If $K>0$ is a constant, denote by $A_{K,h}$ the set of points $x\in B$ for which there exists a bounded linear form $\varphi_x\colon E\to{\mathbb{R}}$ such that, for all $y\in B$,
 \[|h(y)-h(x)-\varphi_x(y-x)|\leq K\|y-x\|^2.\]
 Then the map $h$ has a derivative at each point $x\in A_{K,h}$, and $d_xh=\varphi_x$. Moreover, the restriction of the map $x\mapsto d_xh$ to $\{x\in A_{K,u}\colon \|x\|<1/3\}$ is Lipschitzian with Lipschitz constant $\leq 6K$.
\end{lemma}
As a partial converse we also have
\begin{lemma}\label{lem:lipschitzconverse}
 Let $B$ be the open unit ball in the normed space $E$. If $h\colon B\to {\mathbb{R}}$ is differentiable and the map $x\mapsto dh_x$ is $K$-Lipschitz for some $K>0$, then for all $x,y\in B$, we have
 \[|h(y)-h(x)-dh_x(y-x)|\leq K\|y-x\|^2.\]
\end{lemma}

\begin{lemma}\label{lem:lipcondholds}
 Let $\gamma$ be a curve as in Lemma \ref{lem:existenceLx}, and additionally assume that $\gamma$ is differentiable at 0 in the sense that 0 is a Lebesgue point for $d\gamma\in L^1([-t,t];{\mathbb{R}}^d)$. Let $x=\gamma(0)$.
 Then there is some $K>0$ such that, for $y\in V$,
 \begin{equation}\label{eq:lipcondholds}
  |f(y)-f(x)-\frac{\partial \tilde L}{\partial v}(x,\dot\gamma(0))(y-x)|\leq K\|y-x\|^2
 \end{equation}

\end{lemma}
\begin{proof}
 Let $0<\varepsilon\ll t$.
 Let $D$ be an open subset of $TV$ that contains $\operatorname{supp}\mu\cap TV$ as well as all vectors of size $\leq 1/\varepsilon$, and has compact closure $\overline D$, and let $K_0$ be the Lipschitz constant of $(x,v)\mapsto d\tilde L(x,v)$ on $D$. 
 For  $q\in{\mathbb{R}}^d$, $\|q\|\leq 1$, let $\gamma_{q,\varepsilon}\colon [-\varepsilon,0]\to V$ be the curve such that 
 \[\gamma_{q,\varepsilon}(s)-\gamma(s)=\frac{\varepsilon+s}{\varepsilon}q, \quad s\in[-\varepsilon,0]\]
 and $\gamma_{q,\varepsilon}(s)=\gamma(s)$ for $s\in[-t,-\varepsilon]$.
 By Lemmas \ref{lem:Lsuperdiff}, \ref{lem:existenceLx}, and \ref{lem:lipschitzconverse}, for $q$ close enough to the origin, we have, for almost all $s\in[-t,0]$,
 \begin{equation}\label{eq:lipschitzineq}
  \tilde L(d\gamma_{q,\varepsilon}(s))-\tilde L(d\gamma(s))-\tilde L_x(d\gamma(s))\frac{\varepsilon+s}{\varepsilon}q-\tilde L_v(d\gamma(s))\frac q\varepsilon\leq K_0\left\|\frac q\varepsilon\right\|^2.
 \end{equation}

 Since $\tilde L=\alpha+\tilde g$ with $\tilde g\geq 0$, we have
 \[f(x+q)-f(\gamma(-t))=\int_{-t}^0\alpha_{\gamma_{q,\varepsilon}(s)}(\dot\gamma_{q,\varepsilon}(s))ds\leq \int_{-t}^0 \tilde L(d\gamma_{q,\varepsilon}(s))\,ds\]
 and also, since $f$ satisfies \eqref{eq:laxoleinikpast},%
 \[f(x)-f(\gamma(-t))=\int_{-t}^0 \tilde L(d\gamma(s))\,ds.\]
 Thus,
 \begin{multline}\label{eq:taylorexpansion}
  f(x+q)-f(x)\leq \int_{-t}^0\tilde L(d\gamma_{q,\varepsilon})-\tilde L(d\gamma(s))\,ds\\
  \leq \int_{-\varepsilon}^0\tilde L_x(d\gamma(s))\frac{\varepsilon+s}{\varepsilon}q+\tilde L_v(d\gamma(s))\frac q\varepsilon\,ds+K_0\left\|\frac q\varepsilon\right\|^2,
 \end{multline}
 where the last inequality follows from \eqref{eq:lipschitzineq}. Note that it follows from Lemma \ref{lem:existenceLx} that $\tilde L_x(d\gamma(s))$ is well defined for almost every $s$.
 
 Now we will show that 
 \begin{equation}\label{eq:auxiliarystatement}
  \int_{-\varepsilon}^0\tilde L_x(d\gamma(s))\frac{\varepsilon+s}{\varepsilon}q+\tilde L_v(d\gamma(s))\frac q\varepsilon\,ds=\tilde L_v(d\gamma(0))q.
 \end{equation}

 We let $\psi\colon{\mathbb{R}}\to {\mathbb{R}}$ be a smooth, nonnegative function, vanishing in a neighborhood of 0, and equal to 1 outside a slightly larger neighborhood of 0. We write, for $0<r\ll 1$,
 \begin{align}
  \notag\int_{-\varepsilon}^0&\tilde L_x(d\gamma(s))\frac{\varepsilon+s}{\varepsilon}q+\tilde L_v(d\gamma(s))\frac q\varepsilon\,ds=\\
  \label{eq:termthatgoesaway}&\int_{-\varepsilon}^0\tilde L_x(d\gamma(s))\psi(\tfrac s r)\frac{\varepsilon+s}{\varepsilon}q+\tilde L_v(d\gamma(s))\frac{\partial}{\partial s}[\psi(\tfrac s r)\frac{\varepsilon+s}{\varepsilon}q]\,ds \\
  \label{eq:termthatstays}&+\int_{-\varepsilon}^0\tilde L_x(d\gamma(s))(1-\psi(\tfrac s r))\frac{\varepsilon+s}{\varepsilon}q+\tilde L_v(d\gamma(s))\frac{\partial}{\partial s}[(1-\psi(\tfrac s r))\frac{\varepsilon+s}{\varepsilon}q]\,ds
  \end{align}
  The term \eqref{eq:termthatgoesaway} vanishes because $h(s)=\psi(\tfrac s r)\frac{\varepsilon+s}\varepsilon q$ satisfies the condition for \eqref{eq:nullonvanishing} to hold (it vanishes at $s=-\varepsilon,0$). Now, as $r\to0$, we see that the first term in \eqref{eq:termthatstays} vanishes asymptotically because $1-\psi(\tfrac s r)$ tends to 0. So we are left with the second term in \eqref{eq:termthatstays}, which we expand to get
  \begin{equation}\label{eq:termexpansion}
   -\int_{-\varepsilon}^0\tfrac 1 r \psi'(\tfrac s r)\tilde L_v(d\gamma(s))\frac{\varepsilon+s}{\varepsilon}q\,ds+\int_{-\varepsilon}^0 \tilde L_v(d\gamma(s))(1-\psi(\tfrac s r))\frac{q}{\varepsilon}\,ds.
  \end{equation}
  Now, the second term in \eqref{eq:termexpansion} vanishes as $r\to0$ because, again, $1-\psi(\tfrac s r)$ tends to 0. The first term, on the other hand, contains $-\tfrac1r\psi'(\tfrac s r)$, which approximates the Dirac delta function at $s=0$ as $r\to 0$, so the integral converges to $\tilde L_v(d\gamma(0))q$, which proves \eqref{eq:auxiliarystatement}.
  
  Hence, we have from \eqref{eq:taylorexpansion} and \eqref{eq:auxiliarystatement},
  \[f(x+q)-f(x) \leq \tilde L_v(d\gamma(0))q+\frac {K_0} {\varepsilon^2}\|q\|^2.\]
  A similar argument applied to $\gamma$ on $[0,t]$ and using \eqref{eq:laxoleinikfuture} instead of \eqref{eq:laxoleinikpast}, gives
  \[f(x+q)-f(x) \geq \tilde L_v(d\gamma(0))q-\frac {K_0} {\varepsilon^2}\|q\|^2,\]
  thus proving the lemma with $K=K_0/\varepsilon^2$.
\end{proof}

\subsubsection{Time-dependent setting}
\label{sec:timedependentlipschitz}
Recall that the time dependent setting was defined in Section \ref{sec:timedependent}. The following is a corollary of Theorem \ref{thm:mather}, with $M=N\times P$.
\begin{corollary}\label{cor:mathertimedependent}
 Let $N\times P$ be a manifold that is the product of two $C^\infty$ manifolds $N$ and $P$, with $P$ playing the role of time, $\dim P=1$, so that we are in the time-dependent setting. Let $c$ be a normal 0-current on $N\times P$ such that $\mathscr H_1(c)$ is not empty. Assume that $L$ is an element of $C^2(TN\times P)$ such that $\nu\mapsto\int_{TN\times P}L\,d\nu$ reaches its minimum within $\mathscr H_1(c)$ at some point $\mu$. Let $U$ be an open subset of $N\times P$ with compact closure such that $\pi(\operatorname{supp}\mu)\subseteq U$. Then there exist a Lipschitz function $f\colon U\subseteq N\times P\to{\mathbb{R}}$, a nonnegative function $g\colon TU\cap (TN\times P\times\{\mathbf 1\})\to {\mathbb{R}}_{\geq 0}$, and a bounded (possibly discontinuous) section $\alpha\colon U\subseteq N\times P\to T^*N\times {\mathbb{R}}$ such that:
 \begin{enumerate}
  \item $\alpha$ is a Clarke differential of $f$ (coinciding with $\alpha=df$ wherever $f$ is differentiable;
  \item throughout $TU\cap (TN\times P\times \{\mathbf 1\})$ we have
   \[L=\begin{cases}\int_{TN\times P}L\,d\mu+\alpha+g,&c=0,\\\alpha+g,&c\neq0,\end{cases}\]
  which in particular means that
   \[L=\begin{cases}\int_{TN\times P}L\,d\mu+df+g,&c=0,\\ df+g,&c\neq0,\end{cases}\]
 wherever $f$ is differentiable;
  \item for every open set $V\subset N\times P$ not intersecting the support of $c$, $f$ is differentiable on $\pi(\operatorname{supp}\mu)\cap V$ and, for $(x_0,t_0,v_0,\mathbf 1)\in\operatorname{supp}\mu\cap TV$, we have
   \begin{multline*}
    df_{(x_0,t_0)}(v,\tau)=\frac{\partial L}{\partial v}(x_0,t_0,v_0)\cdot v+\\\left(L(x_0,t_0,v_0)-\frac{\partial L}{\partial v}(x_0,t_0,v_0)\cdot v_0-\int_{TN\times P}L\,d\mu\right)\cdot \tau, \quad c=0,
   \end{multline*}
   and
      \begin{multline*}
    df_{(x_0,t_0)}(v,\tau)=\frac{\partial L}{\partial v}(x_0,t_0,v_0)\cdot v+\\
    \left(L(x_0,t_0,v_0)-\frac{\partial L}{\partial v}(x_0,t_0,v_0)\cdot v_0\right)\cdot \tau, \quad c\neq 0,
   \end{multline*}
   where $v\in T_{x_0}N$, $\tau\in TP$, and $\partial L/\partial v$ denotes the derivative of the restriction of $L$ to $T_{x_0}N\times \{t_0\}\times \{1\}$;
  \item for every open neighborhood $Y\subset U$ of $\operatorname{supp} c$, there is a constant $C_Y>0$ such that the map $(x,t)\mapsto df_{(x,t)}$ is $C_Y$-Lipschitz throughout $\pi(\operatorname{supp}\mu)\cap U\setminus Y$;
  \item $g\equiv 0$ throughout $\operatorname{supp}\mu$.
 \end{enumerate}
 Conversely, if $L$ has this structure, its action reaches its minimum within $\mathscr H_1(c)$ at $\mu$.
\end{corollary}

\section{Examples}
\label{sec:examples}

\begin{example}[Exact form]\label{ex:exactform}
If $L$ is itself an exact form, that is, if $L=df$ for some $f\in C^\infty(M)$, then its action can be minimized by any closed measure $\mu\in \mathscr H(0)$. This shows that it would be impossible to prove any regularity for the minimizers without stronger hypotheses on $L$. It also shows that \emph{every measure $\mu\in \mathscr H(0)$ is a minimizer} of infinitely many Lagrangians, so that it would be hopeless to try to prove the regularity of the minimizers.
\end{example}

\begin{example}[Tonelli Lagrangians]\label{ex:tonelli}
 In the time-dependent setting on $N\times P$ with $P=S^1={\mathbb{R}}/{\mathbb{Z}}$, if $L$ is strictly convex and super-linear in the fibers of $TN$, the existence of minimizers was proved by Tonelli; see for example \cite{fathibook}. From Corollary \ref{cor:mathertimedependent} we recover Mather's theory \cite{matheractionminimizing91} in slightly greater generality because we do not require the minimizers to be invariant under the Euler-Lagrange flow (and we are not the first ones to achieve this greater degree of generality; see \cite{fathibook}): for minimizers in $\mathscr H_1(0)$, the Lipschitzity of $df$ implies in this case that $\operatorname{supp}\mu$ defines a Lipschitz fibration. 
 
 This is a context that has been studied very extensively. Among other results, we mention that it has been proved that $f$ can be chosen to be $C^{1,1}$ throughout $N\times P$ \cite{bernard2007existence} or as a so-called viscosity solution of the Hamilton-Jacobi equation \cite{fathibook}, a property that implies much stronger regularity than we prove in Theorem \ref{thm:mather} and has interesting consequences regarding the associated dynamical system. A good summary can be found in \cite{fathiicm}.
\end{example}

\begin{example}[Irregularity outside $\pi(\operatorname{supp}\mu)$]\label{ex:discontinuousdf}
 Let $M=S^1={\mathbb{R}}/{\mathbb{Z}}$, and let $f\colon S^1\to{\mathbb{R}}$ be a Lipschitz function, differentiable at 0, and with $f'(0)=0$. We let $\mu=\delta_{(0,0)}$. Let $L\colon TS^1\to{\mathbb{R}}$ be a smooth function with $L\geq df$, $L(0,0)=0$, and such that $\int_{S^1} L(x,r)-df_x(r)\,dx\to 0$ as $r\to \pm\infty$. Then this $f$ is the only possible such function in the statement of Theorem \ref{thm:mather}. In the theorem, $f$ is shown to be $C^{1,1}$ on $\pi(\operatorname{supp}\mu)$, and this example shows that no better result can be obtained outside $\pi(\operatorname{supp}\mu)$.
\end{example}

\begin{example}[Sharpness of the Lipschitzity of $df$]\label{ex:sharpregularity}
 This example was communicated to the author by Stefan Suhr, who learned it from Victor Bangert.
 
 In the Beltrami--Klein model of 2-dimensional hyperbolic space, the geodesics correspond to the straight lines on the unit disc $\mathbb D$. Let $g$ be the corresponding Riemannian metric on $\mathbb D$.
 Take the family $\Gamma$ of straight rays that emanate from ${\mathbb{R}}_{\leq 0}\subset {\mathbb{C}}$ and are vertical within $\{{\mathbb{R}}e z\leq 0\}$ and radial from the origin within $\{{\mathbb{R}}e z\geq 0\}$.
 The family $\Gamma$ foliates $\mathbb D\setminus {\mathbb{R}}_{\leq 0}$. Note that the derivatives of the geodesics in $\Gamma$ are only Lipschitz-varying, as the rate of change of these derivatives is not differentiable at $\{{\mathbb{R}}e z=0\}$.  
 
 We consider the case in which $L(x,v)=g_x(v,v)=|v|^2$. Take any closed subset $A$ of $\mathbb D\setminus {\mathbb{R}}_{\leq 0}$ that is bounded in the metric $g$ (in particular, it is bounded away from the circle $\partial \mathbb D$). The geodesics in $\Gamma$ can be indexed by ${\mathbb{R}}$, so that $\Gamma=\{\gamma_r\}_{r\in{\mathbb{R}}}$. Take a measure $\nu$ on ${\mathbb{R}}$. Assuming that we have a unit-speed parameterization of the geodesics $\gamma\in \Gamma$, $g(\dot\gamma,\dot\gamma)=1$, the measure $\mu$ defined by
 \[\int_{T\mathbb D}\phi\,d\mu=\int_{\mathbb{R}}\int_{\{t:\gamma_r(t)\in A\}}\phi(\gamma_r(t),\dot\gamma_r(t))\,dt\,dr,\quad\phi\in C^0(T\mathbb D),\]
 is a minimizer of the action of $g$ within the set of measures that share its boundary, that is, within $\mathscr H(\partial \mu)$. The function $f$ in this case corresponds to the distance to ${\mathbb{R}}_{\leq 0}$.
 
 Since $g$ is smooth, and since necessarily we have 
 \[
  df_xv=g(\dot\gamma_r(t),v)
 \]
 (for $x\in \pi(\operatorname{supp}\mu)$ and for $r$ and $t$ such that $x=\gamma_r(t)$), it follows that $df$ only has Lipschitz regularity in this example. This shows that the version of Theorem \ref{thm:mather} for measures with boundary $\partial\mu$ cannot be improved, and suggests that the same is true for closed measures.
\end{example}

\begin{example}[More irregularity in the wild]\label{ex:arnaud}
 This example was suggested to the author by Marie-Claude Arnaud.
 
 In \cite[\S 4.2]{contrerasfigallirifford}, a Riemannian metric is constructed on the 2-dimensional torus ${\mathbb{T}}^2={\mathbb{R}}^2/{\mathbb{Z}}^2$ that is hyperbolic except inside a small disc $D\subset {\mathbb{T}}^2$. A method is then described to find a measure $\mu$ that minimizes the action and does not correspond to a closed orbit because it has irrational homology. The measure $\mu$ is invariant under the geodesic flow of ${\mathbb{T}}^2$.
 
 The theory developed in \cite{arnaud2011link} can be adapted to analyze the regularity of $\operatorname{supp}\mu$. That theory is about maps on the annulus $S^1\times{\mathbb{R}}$. To adapt it, take a smooth circle $\beta\subset M$ transversal to restriction of the geodesic flow determined by $\operatorname{supp}\mu$, and look at the map $\phi\colon T\beta\to T\beta$ determined by the first-return map of that flow.
 
 What the theory of \cite{arnaud2011link} tells us is that $x\mapsto df_x$ in this case cannot be the the restriction of a $C^1$ section of $T^*{\mathbb{T}}^2$. Already from Mather's theory \cite{matheractionminimizing91} (or from Theorem \ref{thm:mather}) we know that it must be Lipschitz, and the question remains as to whether it is something in-between.
\end{example}
\begin{example}[Young Superposition Principle] \label{ex:youngsuperposition}
 Let $c$ be a normal 0-current on $M$ and write $c=c^+-c^-$ with positive Radon measures $c^\pm$, and assume $\langle c,1\rangle=c^+(M)-c^-(M)=0$, so that $c$ is a boundary.
 
 Let $\mathcal G(c)$ be the set of \emph{generalized curves connecting $\operatorname{supp} c^-$ and $\operatorname{supp} c^+$}, namely, the set of compactly-supported, positive, Radon measures $\mu$ on $TM$ such that there are some $T>0$ and a Lipschitz curve $\gamma\colon[0,T]\to M$ such that $\gamma(0)\in\operatorname{supp} c^-$, $\gamma(T)\in\operatorname{supp} c^+$, and for all $x\in \pi(\operatorname{supp}\mu)$ there is $t\in[0,T]$ with $\gamma(t)=x$ and for almost every $t\in[0,T]$, $\int_{T_xM}v\,d\mu_x=\gamma'(t)$, where $\mu=\int_M \mu_x\,d(\pi_*\mu)$ is the canonical fiberwise disintegration of $\mu$. We take $\mathcal G(0)$ to be the set of measures of the form just described with $\gamma(0)=\gamma(T)$.
 
 Let $\mathcal C(c)$ be the set of measures $\mu_\gamma$ induced by $C^\infty$ curves $\gamma\colon[0,T]\to{\mathbb{R}}$ starting in $\gamma(0)\in\operatorname{supp} c^-$ and ending in $\gamma(T)\in\operatorname{supp} c^+$ by pushing forward the Lebesgue measure ${\operatorname{Leb}}_{[0,T]}$ into $TM$ with the derivative of $\gamma$, that is, $\mu_\gamma=(\gamma,\gamma')_*{\operatorname{Leb}}_{[0,T]}$.  Again, we take the set $\mathcal C(0)$ to be the set of measures of the form just described with $\gamma(0)=\gamma(T)$.

 Here we show how to use Theorem \ref{thm:characterization} to show 
 \begin{proposition}\label{prop:superposition}
  The set $\mathscr H(c)$ from Definition \ref{defn:currents} is the subset of measures $\mu$ satisfying $\partial\mu=c$ in the the convex envelope of the positive cone over $\mathcal G(c)\cup\mathcal G(0)$:
 \[\mathscr H(c)=\{\mu:\partial\mu=c\}\cap\overline{\operatorname{conv}{\mathbb{R}}_{\geq0} (\mathcal G(c)\cup \mathcal G(0))}.\]
 Also, $\mathcal H(c)$ is the subset of measures $\mu$ satisfying $\partial\mu=c$ in the closed convex envelope of the set of positive multiples of measures induced by $C^\infty$ curves $\gamma\colon[0,T]\to{\mathbb{R}}$ that are either closed or satisfy $\gamma(0)\in\operatorname{supp} c^-$ and $\gamma(T)\in\operatorname{supp} c^+$:
 \begin{equation} \label{eq:superposition}
  \mathscr H(c)=\{\mu:\partial\mu=c\}\cap\overline{\operatorname{conv}{\mathbb{R}}_{\geq0}(\mathcal C(c)\cup\mathcal C(0))}.
 \end{equation} 
 \end{proposition}

 This is known as the \emph{Young Superposition Principle} \cite[\S 5]{patrick} since it first appeared in Young's book \cite{young}.
 
 \begin{proof} 
  To prove these statements, observe that the set $\mathcal G(c)$ is closed, and that $\mathcal C(c)$ is dense in it, so the second statement will imply the first one. To prove the second statement, we will use the fact that a closed convex cone $K$ in the set of compactly supported distributions $\mathscr E'$ is completely determined by the set $K'$ of continuous linear functionals that are nonnegative on it; see Remark \ref{rmk:duality}.
  Let $\mathcal A$ be the set of positive, compactly supported, Radon measures $\mu$ on $TM$ whose boundary is a signed measure $\nu$ on $M$, $\nu=\nu^+-\nu^-$, with $\operatorname{supp}\nu^\pm\subset\operatorname{supp} c^\pm$.
  By Theorem \ref{thm:characterization} we know that $\mathcal A$ is the subset of the set of compactly-supported distributions $\nu\in\mathcal E'$ with $\langle \nu,L\rangle\geq 0$ for all $L\in C^\infty(TM)$ of the form
  \[L=C+\lim_{i\to+\infty}df_i+g_i\]
  with $C\geq 0$ and functions $f_i\in C^\infty(M)$ and $g_i\in C^\infty(TM)$ such that 
  \[\lim_{i\to+\infty}\left(\min_{\operatorname{supp} c^+} f_i-\max_{\operatorname{supp} c^-} f_i\right)\geq 0\quad \textrm{and}\quad g_i\geq 0.\] %
  Similarly, Theorem \ref{thm:characterization} shows that the same functionals that are nonnegative on the set $\overline{\operatorname{conv}{\mathbb{R}}_{\geq0}(\mathcal C(c)\cup\mathcal C(0))}$ are exactly the same ones. 
  Observe that the set $\mathscr H(c)$ is the (closed and convex) subset of $\mathcal A$ consisting of measures $\mu$ with $\partial \mu=c$. Taking the intersection in the right-hand side of \eqref{eq:superposition}, we get the equality.
 \end{proof}

\end{example}

\begin{example}[Non-strictly convex Finsler metrics]\label{ex:finsler}
 Denote by $TM_{\neq 0}$ the set of points $(x,v)\in TM$ with $v\neq 0$.
 A \emph{non-strictly convex Finsler metric} is a function $m\colon TM\to{\mathbb{R}}$ that is homogeneous of degree 1 in the fibers, meaning that $m(x,av)=am(x,v)$ for all $a>0$, is everywhere positive on $TM_{\neq 0}$, and is convex on each fiber $T_xM$, $x\in M$. Assume that the manifold $M$ is connected, and that $m$ is $C^{2}$ on $TM_{\neq 0}$. Let $X$ be a closed subset of $M$, and let the distance from $X$ to a point $x\in M$ be defined by
 \[\operatorname{dist}_m(X,x)=\inf_{\gamma}\int_0^Tm(d\gamma(t))dt,\]
 where the infimum is taken over all absolutely-continuous curves $\gamma\colon[0,T]\to M$ with $\gamma(0)\in X$ and $\gamma(T)=x$. 
 The distance $\operatorname{dist}_m(X,x)$ is always realized by a (non-unique) absolutely continuous curve $\gamma_x$. Let $A\subset M\setminus X$ be the set of points $x$ such that there is only one point in $X$ that is at distance $\operatorname{dist}_m(X,x)$ from $x$. All measures that are convex combinations of the measures on $TM$ induced  by integration over the derivatives $d\gamma_x$, $x\in A$, are minimizers in $\mathscr H(\partial \mu)$ of the action of $m$.
 
 For each $m$ of this kind, there exists a Lagrangian $L\colon TM\to{\mathbb{R}}$ that is $C^2$, convex, super-quadratic on the fibers, $L\geq m$, and $L(x,av)=m(x,av)$ for exactly one $a\neq 0$ for each $(x,v)$, $v\neq 0$. Applying Theorem \ref{thm:mather} to this function $L$ together with the minimizers $\mu$, we obtain in fact that the function $f(x)=\operatorname{dist}_m(X,x)$ is of class $C^{1,1}_{\operatorname{loc}}$ throughout $A$.

 This extends some of the results of \cite{linirenberg,linirenberglong} to the non-strictly convex case.
\end{example}

\begin{example}[Finsler-like optimal mass transport]\label{ex:finslertransport}
 Let $m\colon TM\to{\mathbb{R}}$ be a function that is homogeneous of degree 1 in the fibers, meaning that $m(x,av)=a\,m(x,v)$ for all $a>0$, and assume that $m$ is $C^2$ away from the zero section. Let $\nu_1$ and $\nu_2$ be two positive Radon measures on $M$, and consider the following optimal mass transport problem: minimize the cost
 \[\int_{TM} m(x,v)\,d\mu\]
 among all measures $\mu$ on $TM$ whose induced current has boundary $\nu_2-\nu_1$.
 It follows from Theorem \ref{thm:mather} that, \emph{if there is a measure $\mu$ on $TM$ solving that problem, then $m$ must be of the form
 \[m=df+g\]
 for some Lipschitz function $f\colon M\to{\mathbb{R}}$ that is $C^{1,1}_{\operatorname{loc}}$ on $\pi(\operatorname{supp}\mu)\setminus (\operatorname{supp}\nu_1\cup\operatorname{supp}\nu_2)$ and some measurable, nonnegative function $g\colon TM\to{\mathbb{R}}_{\geq0}$ that vanishes identically on $\operatorname{supp}\mu\setminus\pi^{-1}(\operatorname{supp}\nu_1\cup\operatorname{supp}\nu_2)$}.
 
 To apply Theorem \ref{thm:mather} in this case, first note that, due to the 1-homogeneity of $m$, we may assume that the minimizers are supported on the unit sphere of the fibers of $TM$. Then take the function $L\colon TM\to{\mathbb{R}}$ to be $C^2$, equal to $m$ outside a neighborhood $U$ away from the zero section, and slightly greater than $m$ in $U$.
 
 The argument sketched above only gives us information relating to what happens away from $\operatorname{supp}\nu_1\cup\operatorname{supp}\nu_2$; to overcome this issue, one can use the Young Superposition Principle (Example \ref{ex:youngsuperposition}) to conclude that the same result given by Theorem \ref{thm:mather} must be valid for every subset of the generalized curves that compose the minimizer $\mu$, with the same functions $f$ and $g$, and note that (the support of) the boundary of each generalized curve consists of just two points on $M$.
 
 This generalizes some of the results of \cite{bernardbuffonimonge}. We remark also that the result that, when $m^2$ is a Riemannian metric, the flow encoded by $\mu$ actually goes along a gradient vector field \cite{ambrosiogiglisavare} also follows from the statements sketched above: that vector field is precisely the dual of $df$, namely, $\nabla f$.
\end{example}

\section{Optimal control and the maximum principle}
\label{sec:ocontrol}
To illustrate how the methods we have developed in Section \ref{sec:characterization} can also be used to understand problems of optimal control, we will apply them to a slight relaxation of a problem discussed in \cite[\S\textsc{iii}.3]{bardi}, which is known as a \emph{finite horizon Legendre problem}. Although the theory could be developed in greater generality, we refrain from doing this here in the interest of simplicity.%

The main result of this section is Theorem \ref{thm:matherocontrol}, which is analogous to Theorem \ref{thm:mather}, and its content is linked to the previously-existing theory of optimal control in Remark \ref{rmk:HJB}, where we recover the Maximum Principle and the Hamilton--Jacobi--Bellman equation. Also, Theorem \ref{thm:ocontrolcharact} is a coarse characterization of the minimizable functionals analogous to Theorem \ref{thm:characterization}.

The relaxation we will consider replaces ---in an analogous way as we do in our treatment in Section \ref{sec:characterization}--- curves with measures; this is almost equivalent to the relaxation described in \cite[\S\textsc{iii}.2.5 (pp. 113--116)]{bardi} that entails the introduction of so-called \emph{relaxed} or \emph{chattering} controls. As in the statement of \cite[Corollary \textsc{iii}.2.21]{bardi}, this relaxation of the problem is not very significant because the resulting value function should coincide with the one corresponding to the original finite horizon Legendre problem.

We will try to keep the same notations as in the book, although we will immediately translate to the setting that interests us.

\begin{remark}
 The Lagrangians in equations \eqref{eq:defL} and \eqref{eq:convLocontrol} below draw the connection between the approach taken here and the optimal transport theory described in \cite{patrick,bernardbuffoni}. Thus in the present context it is also possible to think about the minimizer as describing a sort of transport plan, except that the velocities must be obtained through the controls. Still, a version of the Young Superposition Principle (Example \ref{ex:youngsuperposition}) holds in this context, so that the measure $\mu$ can be understood as encoding a set of generalized curves describing the trajectories of the transportation of the component of $\partial\mu$ at $t=0$ to its component at $t=t_0$. Also, the Lipschitz regularity results of Theorem \ref{thm:matherocontrol} generalize those obtained in \cite{bernardbuffoni}.
\end{remark}

\subsection{Setting}

We fix a time interval $I=[0,t_0]$, $t_0>0$. We let $t\colon I\to{\mathbb{R}}$ be a chart, and we will denote by $\mathbf 1$ the vector field tangent to $I$ such that $dt(\mathbf 1)=1$.

We let $A$ be a topological space 
to serve as the set of controls, and  $N>0$ denote the dimension of the space ${\mathbb{R}}^N$ of states. 

We assume that we are given a continuous 
function $f\colon {\mathbb{R}}^N\times A\to{\mathbb{R}}^N$ that gives the dynamics; thus, if we had a curve $y\colon I\to{\mathbb{R}}^N$ in the space of states corresponding to a control $\alpha\colon I\to A$, $y$ would satisfy
\[y'(t)=f(y(t),\alpha(t)),\quad t\in I.\]

However, instead of considering such curves, we will consider the set $\mathscr I$ of compactly-supported, Radon, probability measures $\nu$ on ${\mathbb{R}}^N\times I\times A$ satisfying the condition
\[\int_{{\mathbb{R}}^N\times I\times A}du_{(x,t)}(f(x,a),\mathbf 1)\,d\nu(x,t,a)=0\]
for all $u\in C^\infty({\mathbb{R}}^N\times I)$ vanishing on ${\mathbb{R}}^N\times \partial I={\mathbb{R}}^N\times \{0,t_0\}$. This amounts to requiring $\partial \nu$ to be supported at times $0$ and $t_0$. We remark that in the case of the curve $y$ above, this condition would take the form
\begin{multline*}
 \int_0^{t_0} du_{(y(t),t)}(f(y(t),\alpha(t)),\mathbf 1)\,dt=\int_0^{t_0} du_{(y(t),t)}\left(\tfrac{d}{dt}(y(t),t)\right)dt\\
 =\int_0^{t_0} \tfrac{d}{dt}u(y(t),t)\,dt
 =u(y(t_0),t_0)-u(y(0),0)=0.
\end{multline*}

The cost to minimize is
\[J(\nu)=\int_{{\mathbb{R}}^N\times I\times A}\ell(x,t,a)\,d\nu(x,t,a),\quad \nu\in \mathscr I,\]
where we assume that the running cost $\ell\colon {\mathbb{R}}^N\times I\times A\to{\mathbb{R}}$ smooth in the ${\mathbb{R}}^N$ and $I$ variables $x$ and $t$.%

\subsection{Coarse characterization}

\begin{theorem}\label{thm:ocontrolcharact}
 In the setting just described, assume additionally %
 that the cost $J$ reaches its minimum within $\mathscr I$ at the probability measure $\mu$. Then there exist  sequences of functions $u_i\in C^\infty({\mathbb{R}}^N\times I)$ %
 such that %
 $u_i\equiv 0$ on ${\mathbb{R}}^N\times \partial I$, 
 and
  \begin{equation}\label{eq:elldecomposition}
  \ell\geq \lim_{i\to+\infty}du_i\circ (f,\mathbf 1),%
 \end{equation}
 where the limit is taken in the space of continuous functions $C^0({\mathbb{R}}^N\times I \times A)$ with the topology of uniform convergence on compact sets, %
%
%
%
  and 
 \[
  \lim_{i\to+\infty} \int_{{\mathbb{R}}^N\times I\times A}\ell-du_i\circ(f,\mathbf 1)\,d\mu=0.
 \]
Conversely, if $\ell$ has this structure, $J$ reaches its minimum within $\mathscr I$ at $\mu$.
\end{theorem}

\begin{proof}[Sketch of proof of Theorem \ref{thm:ocontrolcharact}]
 One can prove the theorem by following essentially the same ideas as for Theorem \ref{thm:characterization}. A lemma analogous to Lemma \ref{lem:dual} holds in this setting, with the definitions
 \begin{align*}
  Q=&\{\bar \ell\colon C^0(\mathbb R^N\times I\times A)\to\mathbb{R} \mid \textrm{$\ell$ is affine, its linear part is induced by}\\
   &\quad\textrm{evaluation at an element of $C^0({\mathbb{R}}^N\times I\times A)$, and $\langle\bar \ell,\nu\rangle\geq 0$ for all $\nu\in\mathscr I$}\},\\
  R=&\{\bar \ell\colon C^0({\mathbb{R}}^N\times I \times A)
  \to\mathbb R\mid\textrm{$\bar \ell(\nu)\geq \langle \nu,du\circ (f,\mathbf 1)\rangle $ for some $u\in C^\infty({\mathbb{R}}^N\times I)$}\\
  &\textrm{vanishing on ${\mathbb{R}}^N\times\partial I$ and for all $\nu\in \mathscr I$}.\}
 \end{align*}
 and with the same conclusion that $\overline R=Q$, and then the rest of the argument can be adapted easily.
\end{proof}

\subsection{Partial equivalence with Lagrangian action minimization}
We now transfer the problem of minimizing $J$ to a problem of minimizing a Lagrangian action in a time-dependent setting as in Definition \ref{defn:timedependent}. Assume we are in the setting of Theorem \ref{thm:ocontrolcharact}.
Let $\mathcal D={\mathbb{R}}^N\times I\times f({\mathbb{R}}^N\times A)\subset {\mathbb{R}}^N\times I\times {\mathbb{R}}^N$ so that $\mathcal D\times\{\mathbf 1\}\subset{\mathbb{R}}^N\times I\times {\mathbb{R}}^N\times \{\mathbf 1\}\subset T({\mathbb{R}}^N\times I)$. Let $L\colon \mathcal D\to{\mathbb{R}}$ be the possibly-discontinuous function given by
\begin{equation}\label{eq:defL}
 L(x,t,\bar v)=\inf\{\ell(x,t,a):a\in A,f(x,a)=\bar v\}
\end{equation}
and $\bar \mu$ be the measure on $\mathcal D\times\{\mathbf 1\}$ given by
\[\int_{{\mathbb{R}}^N\times I\times{\mathbb{R}}^N\times \{\mathbf 1\}} \phi\,d\bar\mu=\int_{{\mathbb{R}}^N\times I\times A} \phi(x,t,f(x,a),\mathbf 1)\,d\mu(x,t,a)\]
for measurable functions $\phi\colon T({\mathbb{R}}^N\times I)\to{\mathbb{R}}$. Note that, with these definitions, $\bar\mu$ is the minimizer of the action of $L$,
\[\bar\nu\mapsto\int_{{\mathbb{R}}^N\times I\times f({\mathbb{R}}^N\times A)\times\{\mathbf 1\}}L(x,t,\bar v)\,d\bar\mu(x,t,\bar v,\mathbf1).\]
within the set of measures $\bar\nu$ on $T({\mathbb{R}}^N\times I)$ that are supported on the set $\mathcal D\times\{\mathbf 1\}$ and have boundary $\partial\bar\nu$ contained in ${\mathbb{R}}^N\times\partial I$. 

We observe that the decomposition \eqref{eq:elldecomposition} means that 
\[L(x,t,\bar v)=\lim_{i\to+\infty}(du_{i,x}(\bar v,\mathbf 1)+\bar w_i(x,t,\bar v)), \quad (x,t,\bar v)\in \mathcal D,\]
where
\[\bar w_i(x,t,\bar v)=\inf\{\ell(x,t,a)-du_i(\overline v,\mathbf 1)%
 :a\in A, f(x,a)=\bar v\}\geq 0.\]

\subsection{Lipschitzity}
The partial equivalence between the minimization of $J$ and the minimization of the action of $L$, together with the results we have obtained for the latter and their proof, suggest that in order to obtain results on the regularity of the value function (defined in  \eqref{eq:valuefunction} and discussed in further depth below), assumptions must be made that will ensure first the regularity of the fiberwise convexification $\tilde L$ (defined in \eqref{eq:convLocontrol}) of %
$L$. %

In this direction, we present the following result, whose technical-looking conditions are relatively mild; see Example \ref{ex:mildconditions}.

\begin{theorem}\label{thm:matherocontrol}
 In the setting described above, assume also that the cost $J$ reaches its minimum within $\mathscr I$ at the probability measure $\mu$. Additionally, assume that  
 the convexified Lagrangian function
 \begin{equation}\label{eq:convLocontrol}
  \tilde L(x,t,\bar v)=\sup\{r+\theta(\bar v):r\in{\mathbb{R}},\,\theta\in T_x^*M,\;r+\theta(f(x,a))\leq \ell(x,t,a)\;\forall a\in \mathcal A\}
 \end{equation}
 is finite and $C^{1,1}$ in a neighborhood $W$ of the support of the minimizer $\mu$ and satisfies\footnote{This condition corresponds to the conclusion of Lemma \ref{lem:Lsuperdiff}.} that for every smooth section $\sigma\colon {\mathbb{R}}^N\times I\to{\mathbb{R}}^N$ there are $C>0$ and a smooth form $\eta\in \Omega^1({\mathbb{R}}^N\times I)$ with
  \begin{multline}\label{eq:Ltildeassumption}
   \tilde L(x_1,t_1,\sigma(x_1,t_1))-\tilde L(x_2,t_1,\sigma(x_2,t_2))-\eta(x_1-x_2,t_1-t_2)
   \\\leq C(|x_1-x_2|^2+|t_1-t_2|^2)
  \end{multline}
  for all $(x_1,t_1)$ and $(x_2,t_2)$ such that $(x_1,t_1,\sigma(x_1,t_1))$ and $(x_2,t_2,\sigma(x_2,t_2))$ are contained in $W$.

 Let $U$ be an open subset of ${\mathbb{R}}^N$ with compact closure $\overline U$ such that 
 \[\pi_{{\mathbb{R}}^N\times I}(\operatorname{supp}\mu)\subset U\times I.\] 
 Then there exist a Lipschitz function $u\colon U\times I\to {\mathbb{R}}$, a nonnegative function $w\colon U\times I\times A\to {\mathbb{R}}_{\geq 0}$, and a (possibly discontinuous) bounded section $\beta\colon U\times I\to{\mathbb{T}}^*U$ such that
 \begin{enumerate}
  \item $\beta$ is a Clarke differential of $u$ (in particular, $du=\beta$ whenever $u$ is differentiable);
  \item for all $x\in U, a\in A, t\in I$, we have
   \[\ell(x,t,a)=\beta_{(x,t)}(f(x,a),\mathbf1) +w(x,t,a);\]
  \item $u$ is differentiable on $\pi_{{\mathbb{R}}^N\times I}(\operatorname{supp}\mu)\cap (U\times (0,t_0))$;\footnote{In fact, if $A$ is a subset of ${\mathbb{R}}^n$ and $a\mapsto f(x,a)$ is injective and $C^1$, for $(x,t,a)\in\operatorname{supp}\mu$, $0<t<t_0$, $\bar v=f(x,a)$, $v\in T_xU\cong {\mathbb{R}}^N$, $\tau\in T_tI\cong {\mathbb{R}}$, we have
    \begin{align*}
     &du_{(x,t)}(v,\tau)=\frac{\partial L}{\partial \bar v}(x,t,\bar v)\cdot v+(L(x,t,\bar v)-\frac{\partial L}{\partial \bar v}(x,t,\bar v)\cdot \bar v)\cdot \tau
     \\
     &=\frac{\partial\ell}{\partial a}(x,t,a)\frac{\partial\phi}{\partial \bar v}(x,t,\bar v)\cdot v + (\ell(x,t,a)-\frac{\partial\ell}{\partial a}(x,t,a)\frac{\partial\phi}{\partial \bar v}(x,t,\bar v)\cdot \bar v)\cdot\tau,%
    \end{align*}
    where $\phi(x,t,f(x,a))=a$, %
    so that $\phi_{\bar v}$ is the pseudoinverse of $f_a$. This is analogous to item \ref{it:differentiability} in Corollary \ref{cor:mathertimedependent}.
    }
  \item for $0<a<b<t_0$, on $\pi_{{\mathbb{R}}^N\times I}(\operatorname{supp}\mu)\cap (U\times [a,b])$, the map $(x,t)\mapsto du_{(x,t)}$ is Lipschitz;
  \item \label{it:wvanishes} $w\equiv 0$ throughout $\operatorname{supp}\mu$ (this amounts to the Hamilton--Jacobi--Bell\-man equation together with the Maximum Principle; see Remark \ref{rmk:HJB}).
 \end{enumerate}
Conversely, if $\ell$ has this structure, $J$ reaches its minimum within $\mathscr I$ at $\mu$.

\end{theorem}
\begin{proof}[Sketch of proof]
 The proof of Theorem \ref{thm:matherocontrol} is the same as for Theorem \ref{thm:mather} working with $\tilde L$ and the decomposition from Theorem \ref{thm:ocontrolcharact}; the assumptions on $\tilde L$ were taken in order for this to work.
 
 After replacing $\tilde L$ with $\tilde L$ and using the Arzel\`a--Ascoli theorem to obtain $u$ and $\beta$, and hence also $w=\tilde L-\beta$, we have that $\tilde L|_{T_xU}$ is $C^{1,1}_{\mathrm{loc}}$ for each $x\in U$,  and $(x,v)\mapsto \frac{\partial\tilde L}{\partial v}(x,v)$ is locally Lipschitz by Lemma \ref{lem:Lvlipschitz}. Our assumption \eqref{eq:Ltildeassumption} on $\tilde L$ is equivalent to the conclusion of Lemma \ref{lem:Lsuperdiff}. Lemmas \ref{lem:lipschitzcriterion} and \ref{lem:lipschitzconverse} are very general and do not need to be changed. Lemmas \ref{lem:smirnov}, \ref{lem:existenceLx}, \ref{lem:lipcondholds} and Corollary \ref{coro:laxoleinik} have obvious analogues, and the mechanism of the proof is the same as described for Theorem \ref{thm:mather}.
\end{proof}
\begin{example}\label{ex:mildconditions}
 The following set of conditions imply the hypotheses of Theorem \ref{thm:matherocontrol} and are perhaps simpler to check:
 Assume that $\ell$ is $C^2$, that $A$ is a subset of ${\mathbb{R}}^n$ for some $0\leq n\leq N$, and that $f$ is $C^2$ and satisfies
 \begin{equation*}%
  C(x)|a-b|\leq|f(x,a)-f(x,b)|\leq\frac{1}{C(x)}|a-b|
 \end{equation*}
 for some continuous function $C\colon {\mathbb{R}}^N\to(0,1)$ and for all $a,b\in A$ and all $x\in{\mathbb{R}}^N$. Assume that the cost $J$ reaches its minimum within $\mathscr I$ at the probability measure $\mu$. 
 Additionally, assume that the support of the measure $(\pi_A)_*\mu$ does not intersect the boundary of $A$.
\end{example}

\subsection{Value function}
 Assume that we are in the setting of Theorem \ref{thm:matherocontrol}.
 Recall that the \emph{value function} $v$ is defined \cite[p. 148]{bardi} by 
 \begin{equation}\label{eq:valuefunction}
  v(x,t)=\inf_{\alpha\in \mathcal A}\int_0^{t}\ell(y_x(s,\alpha),s,\alpha(s))ds,\quad (x,t)\in {\mathbb{R}}^N\times I,
 \end{equation}
 where the infimum is taken over the set of measurable functions $\alpha\colon I\to{\mathbb{R}}$ and the curve $y_x$ satisfies $y_x(s,\alpha)=x+\int_0^s f(y_x(r,\alpha),\alpha(r))\,dr$, in other words, $y_x$ is controlled by $\alpha$.
 As usual, the function $v$ satisfies a Dynamic Programming Principle, namely, for $0<\tau<t$,
  \[v(x,t)=\inf_{\alpha\in\mathcal A}\int_{t-\tau}^t \ell(y_x(s,\alpha),s,\alpha(s))ds+v(y_x(\tau,\alpha),t-\tau).\]
 Let $y\colon I\to{\mathbb{R}}^N$ be an absolutely-continuous minimizing curve whose image is contained in $\pi_{{\mathbb{R}}^N\times I}(\operatorname{supp}\mu)$; such a $y$ exists by the results of \cite{smirnov1993decomposition}.
 It follows from \eqref{eq:valuefunction} that %
 the function $v$ satisfies
 \begin{equation}\label{eq:reluv}
  v(y(0),t)=u(y(t),t)-u(y(0),0),
 \end{equation}
 with $u$ as in Theorem \ref{thm:matherocontrol}.

\begin{remark}\label{rmk:HJB}
 Let again $y\colon I\to{\mathbb{R}}^N$ be an absolutely-continuous minimizing curve corresponding to a control $\alpha\in \mathcal A$ such that $(y(t),\alpha(t))\in\operatorname{supp}\mu$ for all $t\in I$. The equation 
 \begin{equation}\label{eq:myHJB}
  L(y(t),t,f(y(t),\alpha(t)))=du_{y(t)}f(y(t),\alpha(t)),
 \end{equation}
 which is true for all $t\in I$ because of item \ref{it:wvanishes} in Theorem \ref{thm:matherocontrol},
 is equivalent to the \emph{Hamilton-Jacobi-Bellman equation}, 
 \begin{equation}\label{eq:originalHJB}
  v_t+H(x,t,d_xv)=0,
 \end{equation}
 with $v$ as in \eqref{eq:valuefunction} and 
 \[H(x,t,p)=\sup_{a\in A}(-f(x,a)\cdot p-\ell(x,t,a)).\]
 This follows from the fact that, for almost all $t\in I$,
 \begin{align*}
  du_{(y(t),t)}&(f(y(t),\alpha(t)),\mathbf 1) =\frac{d}{dt}[u(y(t),t)-u(y(0),0)],\\
  &=\frac{d}{dt}v(y(0),t),\qquad\textrm{by \eqref{eq:reluv}},\\
  &=\frac{d}{dt}[v(y(0),t_0)-v(y(t),t_0-t)],\qquad\textrm{by \eqref{eq:valuefunction}},\\
  &=-v_x(y(t),t_0-t)y'(t)+v_t(y(t),t_0-t)\\
  &=v_t+H(x,t_0-t,d_xv)+\ell
 \end{align*}
 together with the identity \eqref{eq:defL}, which here amounts to 
 \[L(y(t),t,f(y(t),\alpha(t)))=\ell(y(t),t,\alpha(t)).\]
 
 If we do not restrict to the support of the minimizer $\mu$, equation \eqref{eq:originalHJB} becomes an inequality that corresponds to the fact that $w\geq 0$ and can be written
 \[H(x,f(x,a))\leq -v_t=H(x,d_xv),\quad a\in A.\]
 Thus, we see that $H$ reaches its maximum at the points in the support of $\mu$, where this inequality must be the equality \eqref{eq:originalHJB}. This phenomenon is known as the \emph{Pon\-tryag\-in--Boltyanskii--Gamkrelidze--Mishchenko Maximum Principle} after \cite{pontryagin1962}. We observe that, in this setting, the maximum principle is true for all $t\in I$, rather than only for almost-every $t$, as it is usually formulated.
\end{remark}

\section*{Acknowledgements}
I am deeply grateful to Patrick Bernard for his support, for his patience in listening to sometimes very confused presentations of these results, and for his help in clarifying my ideas with numerous questions and suggestions. I am deeply grateful to Marie-Claude Arnaud, Victor Bangert, Jaime Bustillo, Albert Fathi, Uwe Helmke, and Stephan Suhr for their numerous suggestions and discussions. %
I am very grateful to the \'Ecole Normale Superieure de Paris and the Universit\'e de Paris -- Dauphine for their hospitality and support for development of this research. Finally, I am very grateful to the anonymous referee for suggesting interesting connections to parts of the theory that I was unaware of, and for providing other helpful comments.
\section*{List of changes}\label{sec:changes}
\begin{itemize}
 \item The main issue this correction is addressing was an error in the proof of Theorem \ref{thm:characterization}. The proof is now corrected, and for this the statement and proof of Lemma \ref{lem:dual} has been adapted considerably, and the proof of Theorem \ref{thm:characterization} has also changed. The main idea, however, remains the same.
 \item The statements of Theorems \ref{thm:characterization}, \ref{thm:mather}, \ref{thm:ocontrolcharact} and \ref{thm:matherocontrol} and Corollaries \ref{coro:timedepcoarse} and \ref{cor:mathertimedependent} have been reworked; we have basically removed $\lim_{i\to+\infty}\langle c,f_i\rangle=0$ and separated the two different cases for $c=0$ and $c\neq 0$. 
 \item The variable $c_0$ that appeared in Section \ref{sec:ocontrol} before has been removed everywhere.
 \item The last displayed equation in Remark \ref{rmk:HJB} was also corrected.
\end{itemize}


 \providecommand{\href}[2]{#2}
\providecommand{\arxiv}[1]{\href{http://arxiv.org/abs/#1}{arXiv:#1}}
\providecommand{\url}[1]{\texttt{#1}}
\providecommand{\urlprefix}{URL }

\medskip
Received xxxx 20xx; revised xxxx 20xx.
\medskip

\end{document}